\documentclass[oneside, reqno] {amsart}

\usepackage{xypic}
\input xy
\xyoption{all}
\usepackage{epsfig}
\usepackage{amsthm}
\usepackage{amssymb}
\usepackage{amsmath}
\usepackage{amscd}
\usepackage{color}
\usepackage[T1]{fontenc}
\usepackage{upgreek}
\usepackage[font=scriptsize]{caption}
\usepackage{wrapfig}
\usepackage{graphicx}
\usepackage{subfigure}

\DeclareMathOperator*{\esssup}{ess\,sup}

\newcommand{\bg}{\begin{equation}}
\newcommand{\ed}{\end{equation}}
\newcommand{\bga}{\begin{eqnarray}}
\newcommand{\eda}{\end{eqnarray}}
\newcommand{\pf}{\textbf{Proof:\ }}

\def\cbdu{\par{\raggedleft$\Box$\par}}

\newtheorem {Theorem}  {Theorem}

\numberwithin{Theorem}{section}

\newtheorem {Lemma}[Theorem]  {Lemma}

\theoremstyle{definition}

\theoremstyle{remark}
\newtheorem{Remark}[Theorem]{\bf Remark}

\expandafter\chardef\csname pre amssym.def
at\endcsname=\the\catcode`\@ \catcode`\@=11
\def\undefine#1{\let#1\undefined}
\def\newsymbol#1#2#3#4#5{\let\next@\relax
 \ifnum#2=\@ne\let\next@\msafam@\else
 \ifnum#2=\tw@\let\next@\msbfam@\fi\fi
 \mathchardef#1="#3\next@#4#5}
\def\mathhexbox@#1#2#3{\relax
 \ifmmode\mathpalette{}{\m@th\mathchar"#1#2#3}%
 \else\leavevmode\hbox{$\m@th\mathchar"#1#2#3$}\fi}
\def\hexnumber@#1{\ifcase#1 0\or 1\or 2\or 3\or 4\or 5\or 6\or 7\or 8\or
 9\or A\or B\or C\or D\or E\or F\fi}

\font\teneufm=eufm10 \font\seveneufm=eufm7 \font\fiveeufm=eufm5
\newfam\eufmfam
\textfont\eufmfam=\teneufm \scriptfont\eufmfam=\seveneufm
\scriptscriptfont\eufmfam=\fiveeufm

\catcode`\@=\csname pre amssym.def at\endcsname

\newcounter{remark}
\def  \12  {{\frac{1}{2}}}

\def\build#1_#2^#3{\mathrel{\mathop{\kern 0pt#1}\limits_{#2}^{#3}}}

\numberwithin{equation}{section}

\begin{document}
%\currannalsline{0}{2006}

\title[EMHD self-similar singularity]{Self-similar singularities for electron MHD}

%\author{hello}

\author [Mimi Dai]{Mimi Dai}

\address{Department of Mathematics, Statistics and Computer Science, University of Illinois at Chicago, Chicago, IL 60607, USA}
\email{mdai@uic.edu} 

\author [Hannah Guerra]{Hannah Guerra}

\address{Department of Mathematics, Statistics and Computer Science, University of Illinois at Chicago, Chicago, IL 60607, USA}
\email{hguer2@uic.edu} 

\author [Chao Wu]{Chao Wu}

\address{Department of Mathematics, Statistics and Computer Science, University of Illinois at Chicago, Chicago, IL 60607, USA}
\email{cwu206@uic.edu}

\thanks{The authors are partially supported by the NSF grants DMS--2009422 and DMS--2308208. }

\begin{abstract}
We study several types of self-similar solutions for the electron magnetohydrodynamics (MHD) without resistivity, including locally self-similar solutions and pseudo-self-similar solutions. We show that under certain conditions, these types of self-similar blowup solutions can be excluded. 

\bigskip


\hspace{0.02cm}CLASSIFICATION CODE: 35Q35, 76B03, 76D09, 76E25, 76W05.
\end{abstract}

\maketitle

\section{Introduction}

\medskip

\subsection{Overview}

We consider the electron magnetohydrodynamics (MHD) 
\begin{equation}\label{emhd}
\begin{split}
B_t+ \nabla\times ((\nabla\times B)\times B)=&\ \nu\Delta B,\\
\nabla\cdot B=&\ 0
\end{split}
\end{equation}
which serves as an approximating model for the full MHD system with Hall effect when the motion of the ion flow is slow and can be neglected (cf. \cite{ADFL, Bis1}). The unknown vector $B$ stands for the magnetic field and $\nu\geq 0$ is the resistivity parameter.  We also denote $J=\nabla\times B$ by the current density. The rapid magnetic reconnection phenomena in plasmas is captured by the Hall effect which results in the presence of the nonlinear term in (\ref{emhd}) (cf. \cite{BDS}). The highly singular Hall term is the source of many interesting yet challenging mathematical problems for the Hall MHD and electron MHD. In particular we note the first equation of \eqref{emhd} is quasi-linear and supercritical, both natures presenting serious barriers in mathematical analysis for nonlinear equations. We limit our effort to the topic of self-similar singularity formation for solutions to \eqref{emhd} in this article. 

The authors of \cite{CWeng} showed that the Hall MHD is either ill-posed, in the sense of norm inflation in some Sobolev space with high regularity, or locally well-posed and the solution develops singularity at a finite time. This result applies to the electron MHD \eqref{emhd} as well. Non-unique weak solutions in Leray-Hopf class for the 3D Hall MHD were constructed in \cite{Dai-non} using a scheme of convex integration. Moreover, 
strong ill-posedness phenomena for the Hall MHD and electron MHD were also discovered in the works \cite{JO-ill2, JO-ill1}. Such results suggest the likelihood of singular behavior of solutions to the Hall/electron MHD due to the Hall term. Nevertheless, it is rather subtle to fully understand the Hall term. As a contrast to the ill-posedness results, we showed \cite{Dai-steady} that the electron MHD with resistivity only in the vertical direction has a global regular solution near the steady state $(0,0,1)$; the authors of \cite{JO-well} proved local well-posedness for the electron MHD without resistivity for large perturbations of nonzero background magnetic fields. 

In this paper, we investigate the possibility of singularity formation in the class of self-similar solutions. Our interest in self-similar solutions for \eqref{emhd} stems from the natural scaling property of the equation. Let $B(x,t)$ be a solution of \eqref{emhd} with the initial data $B_0(x,t)$. In the case $\nu=0$, the rescaled function by a parameter $\lambda$
\begin{equation}\label{scale}
B_\lambda(x,t)=\lambda^\alpha B(\lambda x, \lambda^{\alpha+2}t), \ \ \forall \ \ \alpha\in \mathbb R
\end{equation}
is also a solution of \eqref{emhd} with the initial data $\lambda^\alpha B_0(\lambda x)$. While for $\nu>0$, the scaling is 
\[B_\lambda(x,t)= B(\lambda x, \lambda^{2}t).\]
We focus on the case with zero resistivity, i.e. $\nu=0$.

%The electron MHD has been studied extensively by physicists, mostly through numerical simulations, in order to understand the geometry configurations in the reconnection region of the plasma. For instance, we refer the reader to \cite{RDD, SDS} and references therein. In particular, the steady states of (\ref{emhd}) play a vital role in the investigations \cite{DKM, GH, WHL}. Indeed, equation (\ref{emhd}) has a rich class of equilibria, in both of the 3D and 2D settings. Due to the fact that the 2D geometry is more manageable in numerical study, many contributions \cite{CSZ, KC, WH} in the physics community have been emerged for 2D electron MHD. 
%for instance see \cite{JO, WH}.
%From another point of view, we observe (\ref{emhd}) is quasi-linear and supercritical (see \cite{Dai-W}) which presents obstacles for mathematical analysis of the model. Therefore we focus on the electron MHD in 2D geometry in this article.

The scaling \eqref{scale} suggests that we can consider solutions to \eqref{emhd} with $\nu=0$ in the self-similar form
\begin{equation}\label{self-blowup1}
B(x,t)=\frac{1}{(T-t)^{\frac{\alpha}{\alpha+2}}}H\left( \frac{x}{(T-t)^{\frac{1}{\alpha+2}}}\right), \ \ \alpha>-2
\end{equation}
with the profile vector field $H$ satisfying the equations
\begin{equation}\label{eq-H}
\begin{split}
\frac{\alpha}{\alpha+2}H+\frac{1}{\alpha+2} y\cdot \nabla H+\nabla\times ((\nabla\times H)\times H)=&\ 0,\\
\nabla\cdot H=&\ 0
\end{split}
\end{equation}
where $y=\frac{x}{(T-t)^{\frac{1}{\alpha+2}}}$.
The existence of a non-trivial solution $H$ of \eqref{eq-H} corresponds to the existence of a self-similar solution $B$ of \eqref{emhd} in the form \eqref{self-blowup1} that blows up at time $T$. The self-similar blowup may occur locally near a point $(x_0, T)$, in which case we consider the local self-similar form
\begin{equation}\label{local-blowup}
B(x,t)=\frac{1}{(T-t)^{\frac{\alpha}{\alpha+2}}}H\left( \frac{x-x_0}{(T-t)^{\frac{1}{\alpha+2}}}\right), \ \ (x,t)\in B_{\rho_0}(x_0)\times (t_0, T)
\end{equation}
for some $\rho_0>0$ and $t_0\in(0,T)$. The profile function $H$ satisfies the same equation \eqref{eq-H} with $y=\frac{x-x_0}{(T-t)^{\frac{1}{\alpha+2}}}$.

We will also study a more general type of self-similar solutions as follows,
\begin{equation}\label{pse-self}
B_{\lambda,\mu}(x,t)=\mu(t)H(\lambda(t)x)
\end{equation}
for functions of time $\mu(t)$ and $\lambda(t)$ and a profile vector field $H$.  It is straightforward to compute that $\lambda$, $\mu$ and $H$ satisfy
\begin{equation}\label{eq-H2}
\frac{\mu'}{\lambda^2\mu^2}H+\frac{\lambda'}{\lambda^3\mu}y\cdot\nabla H+\nabla\times ((\nabla\times H)\times H)=0
\end{equation}
with $y=\lambda(t)x$.  Solution in the form \eqref{pse-self} is referred a pseudo-self-similar solution.

\medskip

\subsection{Main results}

The main purpose of the paper is to show that several types of self-similar blowup solutions of \eqref{emhd} with $\nu=0$ can be ruled out under certain conditions. 

Let $X(a,t)$ be the trajectory mapping 
\begin{equation}\notag
\begin{cases}
\frac{\partial X(a,t)}{\partial t}=-J(X(a,t), t),\\
X(a,0)=a
\end{cases}
\end{equation}
and $A(x,t)=:X^{-1}((x,t))$ be the back to label mapping.

Our first result is
\begin{Theorem}\label{thm-1}
Let $\alpha>-2$ and $\nu=0$. Suppose $B\in C([0,T); C^2(\mathbb R^3, \mathbb R^3))$ is a classical solution of \eqref{emhd}. If  
\begin{itemize}
\item [(i)] the trajectory mapping $X(\cdot, t)$ for $t\in(0,T)$ generated by $-J$ is a $C^1$ diffeomorphism; 
\item [(ii)] the profile $H\not\equiv 0$ and there exists $p_0>0$ such that $H\in L^p(\mathbb R^3)$ for all $p\in(0,p_0)$,
\end{itemize}
then there exists no self-similar blowup solution $B$ to \eqref{emhd} in the form \eqref{self-blowup1}. 
\end{Theorem}

\medskip

The following two results concern the non-existence of locally self-similar solutions of \eqref{emhd} in the form \eqref{local-blowup} under conditions on the profile function $H$ depending on the scaling parameter $\alpha$.
\begin{Theorem}\label{thm-less32}
Let $B$ be a local self-similar classical solution to \eqref{emhd} on $[t_0,T)$ in the form \eqref{local-blowup} for some $t_0\in (0,T)$ and $\rho_0>0$.
Assume $H\in L^p\cap C_{loc}^2$ for some $p>2$ and $\nabla\times H\in L^q$ for $q>1$. If $-2<\alpha\leq \frac3p$ or $\frac32<\alpha<\infty$, then $H\equiv 0$. 
%If $\alpha=\infty$ and the self-similar solution \eqref{local-self} is global, then $H\equiv 0$.
\end{Theorem}

\begin{Theorem}\label{thm-32}
Let $\alpha=\frac32$ and $B$ be the local self-similar solution described in Theorem \ref{thm-less32}. Suppose the profile function $H$ belongs to $L^2(\mathbb R^3)\cap C^2_{loc}$. In addition, assume there exists a $\delta>0$ such that $H$ satisfies
\begin{equation}\label{assu3}
|H(y)|\geq c|y|^{-\frac32-\delta}, \ \ \ \mbox{for} \ \ |y|\gg 1,
\end{equation}
\begin{equation}\label{assu4}
|\nabla\times H(y)|\leq C|y|^{1-\delta}, \ \ \ \mbox{for} \ \ |y|\gg 1
\end{equation}
for some constants $c\geq 0$ and $C>0$. Then $H\equiv 0$.
\end{Theorem}

\medskip

Certain types of pseudo-self-similar solutions can also be ruled out for \eqref{emhd}.
Denote the class of functions
\begin{equation}\label{class-A}
\begin{split}
\mathcal A=&\left\{(\lambda(t), \mu(t))\in C^1(-\infty, T): (\lambda,\mu)\neq \left(a_1(T-t)^{-\frac27},  a_2(T-t)^{-\frac37}\right), \right.\\
&\ \ \ \ \ \ \ \ \ \ \ \ \ \ \ \ \ \left. a_1,a_2\in \mathbb R \right\}.
\end{split}
\end{equation}

\begin{Theorem}\label{thm-pse1}
Assume $(\lambda,\mu)\in \mathcal A$. 
There is no pseudo-self-similar solution $B$ to \eqref{emhd} with $\nu=0$ in the form \eqref{pse-self} which is regular on $[0,T)$ and satisfies
\begin{equation}\label{assu1}
\esssup_{0<t<T} \|B(t)\|_{L^2(\mathbb R^3)}<\infty
\end{equation}
and
\begin{equation}\label{assu2}
\lim_{t\to T^-}\int_0^{t}\|\nabla^2B(\tau)\|_{L^\infty}\, d\tau=\infty.
\end{equation}
\end{Theorem}

\begin{Theorem}\label{thm-pse2}
Assume 
\[\lambda(t)=a_1(T-t)^{-\frac27}, \ \ \ \mu(t)=a_2(T-t)^{-\frac37}, \ \ \ a_1,a_2\in \mathbb R.\] 
There is no pseudo-self-similar solution $B$ to \eqref{emhd} with $\nu=0$ in the form \eqref{pse-self} which is regular on $[0,T)$ and satisfies \eqref{assu1}, such that additionally either
\begin{equation}\label{assu3}
\esssup_{0<t<T_0} \|B(t)\|_{L^p(\mathbb R^3)}<\infty, \ \ \ \mbox{for some} \ \ p>2
\end{equation}
or
\begin{equation}\label{assu4}
\lim_{t\to T_0^-}\int_0^t \|\nabla\times J(x,\tau)\|_{L^\infty} \, d\tau<\infty
\end{equation}
holds.
\end{Theorem}

\begin{Remark}
It is an open question that whether the Beale-Kato-Majda (BKM) type of blowup criterion holds for the electron MHD with or without resistivity, although other classical blowup criteria can be shown for the Hall MHD (cf. \cite{Dai-hmhd-reg}). The obstacle of obtaining the BKM blowup criterion relies on the derivative loss in the nonlinear term of the electron MHD. Nevertheless, the results of Theorem \ref{thm-pse1} and Theorem \ref{thm-pse2} indicate that we have the BKM type of blowup criterion for the electron MHD in the class of self-similar solutions.
\end{Remark}

\medskip

\subsection{Previous results on self-similar blowup solutions}
For the 3D Navier-Stokes equation, Leray \cite{Leray} raised the question whether there exists a non-trivial self-similar blowup solution. It was answered in the negative in the work \cite{NRS}. In the inviscid case, that is for the Euler equation, contributions have been made in many works including \cite{BS, Chae-NSE, Chae07, Chae11, CSh, He, Sch, Xue}. In these papers, several types of self-similar solutions were ruled out under certain conditions. The study of self-similar blowup solutions has also been extended to the surface quasi-geostrophic equation, see \cite{BGM, CX}. 

\medskip

\subsection{Notations}
We denote $C$ by a generic constant which may be different from line to line. We use $\lesssim$ as the inequality $\leq$ up to a constant which does not play a role in the estimates.

The rest of the paper is organized as follows. In Section \ref{sec-thm1} we prove Theorem \ref{thm-1}. Section \ref{sec-local} contributes to the study of locally self-similar solutions, in particular, the proof of Theorem \ref{thm-less32} and Theorem \ref{thm-32}.
In Section \ref{sec-pseudo} we consider pseudo-self-similar solutions and prove Theorem \ref{thm-pse1} and Theorem \ref{thm-pse2}.

\bigskip

\section{Proof of Theorem \ref{thm-1}}
\label{sec-thm1}

The conclusion of Theorem \ref{thm-1} is an immediate consequence of the following theorem.

\begin{Theorem}\label{thm-gen}
Suppose $B\in C([0,T); C^2(\mathbb R^3))$ is a classical solution of \eqref{emhd} with $\nu=0$. Assume the trajectory mapping $X(\cdot, t)$ generated by $-J$ for $t\in(0,T)$ is a $C^1$ diffeomorphism. Moreover, assume 
\[B(x,t)=\Psi(t) H(\Phi(t) x), \ \ t\in[0,T)\]
with 
\begin{itemize}
\item [(i)]
$\Psi\in C([0,T); (0,\infty))$, $\Phi\in C([0,T); \mathbb R^{3\times 3})$ and $\det (\Phi(t))\neq 0$ on $[0,T)$;
\item [(ii)] $H\in L^p(\mathbb R^3)$, $\forall p\in (0, p_0)$ for some $p_0>0$.
\end{itemize}
Then we have either $\det(\Phi(t))\equiv \det(\Phi(0))$ on $[0,T)$ or $H\equiv 0$.
\end{Theorem}
\pf
The electron MHD \eqref{emhd} with $\nu=0$ can be written as
\begin{equation}\notag
\partial_t B-J\cdot\nabla B=-B\cdot\nabla J.
\end{equation}
Taking the dot product of the equation with $B$ and rearranging the terms gives 
\begin{equation}\label{eq-norm-B}
\partial_t|B|+(-J\cdot\nabla) |B|=\left( \frac{B}{|B|}\cdot\nabla\right)(-J)\cdot \frac{B}{|B|}|B|.
\end{equation}
Denote 
\begin{equation}\label{S}
S=\left( \frac{B}{|B|}\cdot\nabla\right)(-J)\cdot \frac{B}{|B|}.
\end{equation}
 Considering equation \eqref{eq-norm-B} along the trajectory $X(a,t)$ we have
\[\partial_t|B(X(a,t),t)|=S(X(a,t), t)|B(X(a,t),t)| \]
whose solution is given by
\begin{equation}\label{norm-B}
|B(X(a,t),t)|=|B_0(a)|e^{\int_0^t S(X(a,\tau), \tau)\, d\tau}.
\end{equation}
In view of \eqref{S} we have 
\begin{equation}\label{S-size}
-\|\nabla J(t)\|_{L^\infty}\leq S(x,t)\leq \|\nabla J(t)\|_{L^\infty}, \ \ \ \forall \ \ x\in \mathbb R^3.
\end{equation}
It follows from \eqref{norm-B} and \eqref{S-size} that
\begin{equation}\label{norm-bound}
|B_0(a)|e^{-\int_0^t \|\nabla J(\tau)\|_{L^\infty}\, d\tau}\leq |B(X(a,t),t)|\leq |B_0(a)|e^{\int_0^t \|\nabla J(\tau)\|_{L^\infty}\, d\tau}.
\end{equation}
Applying the back to label mapping to \eqref{norm-bound} yields
\begin{equation}\label{norm-bound-2}
|B_0(A(x,t))|e^{-\int_0^t \|\nabla J(\tau)\|_{L^\infty}\, d\tau}\leq |B(x,t)|\leq |B_0(A(x,t))|e^{\int_0^t \|\nabla J(\tau)\|_{L^\infty}\, d\tau}.
\end{equation}
For the pseudo-self-similar ansatz 
\[B(x,t)=\Psi(t) H(\Phi(t) x),\]
%with some functions of time $\Psi(t)$ and $\Phi(t)$, 
we have
\[B_0(x)=B(x,0)=\Psi(0) H(\Phi(0) x)\]
and hence 
\[H(x)=\Psi(0)^{-1}B_0([\Phi(0)]^{-1}x).\]
Therefore we can rewrite $B(x,t)$ as 
\[B(x,t)=\Psi(t)\Psi(0)^{-1}B_0([\Phi(0)]^{-1}\Phi(t)x).\]
Denote $G(t)=\Psi(t)\Psi(0)^{-1}$ and $F(t)=[\Phi(0)]^{-1}\Phi(t)$. Applying \eqref{norm-bound-2} for such pseudo-self-similar form of $B(x,t)$ we obtain
\begin{equation}\notag
|B_0(A(x,t))|e^{-\int_0^t \|\nabla J(\tau)\|_{L^\infty}\, d\tau}\leq G(t)|B_0(F(t)x)|\leq |B_0(A(x,t))|e^{\int_0^t \|\nabla J(\tau)\|_{L^\infty}\, d\tau}.
\end{equation}
Taking the $L^p$ norm on the last inequality gives
\begin{equation}\notag
\|B_0\|_{L^p}e^{-\int_0^t \|\nabla J(\tau)\|_{L^\infty}\, d\tau}\leq G(t)\left[\det(F(t)) \right]^{-\frac1p}\|B_0\|_{L^p}\leq \|B_0\|_{L^p}e^{\int_0^t \|\nabla J(\tau)\|_{L^\infty}\, d\tau}.
\end{equation}
If $H\not\equiv 0$, then $B_0\not\equiv 0$ and hence 
\begin{equation}\notag
e^{-\int_0^t \|\nabla J(\tau)\|_{L^\infty}\, d\tau}\leq G(t)\left[\det(F(t)) \right]^{-\frac1p}\leq e^{\int_0^t \|\nabla J(\tau)\|_{L^\infty}\, d\tau}.
\end{equation}
Suppose there exists a time $t_1\in(0,T)$ such that $\det(F(t_1))\neq 1$. Taking the limit $p\to 0$, it implies 
\[\int_0^t\|\nabla J(\tau)\|_{L^\infty}\,d \tau=\infty\]
which is a contradiction with the assumption $B\in C([0,T); C^2(\mathbb R^3))$. It completes the proof.

\cbdu

\textbf{Proof of Theorem \ref{thm-1}:}
Applying Theorem \ref{thm-gen} with 
\[\Phi(t)=\frac1{(T-t)^{\frac1{\alpha+2}}}I, \ \ \ \Psi(t)=\frac1{(T-t)^{\frac{\alpha}{\alpha+2}}}\]
where $I$ is the unit matrix in $\mathbb R^{3\times 3}$. For $\alpha>-2$, it is easy to see
\[\det(\Phi(0))=T^{-\frac3{\alpha+2}}.\] 
However, for $t\in (0,T)$, we have 
\[\det(\Phi(t))=(T-t)^{-\frac3{\alpha+2}}\neq \det(\Phi(0)).\]
It thus follows from Theorem \ref{thm-gen} that $H\equiv 0$. Hence there is no self-similar blowup solution $B$ to \eqref{emhd} in the form \eqref{self-blowup1}.

\bigskip

\section{Local self-similar solutions}
\label{sec-local}

In this section we consider the locally self-similar solution in the form \eqref{local-blowup} and provide a proof of Theorem \ref{thm-less32} and Theorem \ref{thm-32}.
Without loss of generality, we fix $x_0=0$ in \eqref{local-blowup}.

\subsection{Local energy inequality}
Assume $B$ is regular such that the local energy equality on the region of self-similar is satisfied, i.e.
\begin{equation}\label{en-local1}
\int_{t_1}^{t_2}\int_{\mathbb R^3} \partial_tB\cdot B\sigma \, dxdt+\int_{t_1}^{t_2}\int_{\mathbb R^3} \nabla\times ((\nabla\times B)\times B)\cdot B\sigma \, dxdt=0
\end{equation}
for any $0<t_1<t_2<T$ and the test function $\sigma\in C_0^\infty((0,T)\times \mathbb R^3)$. Note the first integral in \eqref{en-local1} can be written as
\begin{equation}\notag
\begin{split}
\frac12\int_{t_1}^{t_2}\int_{\mathbb R^3} \partial_t|B|^2\sigma \, dxdt=&\ \frac12\int_{\mathbb R^3} |B(x,t_2)|^2\sigma(x,t_2) \, dx-\frac12\int_{\mathbb R^3} |B(x,t_1)|^2\sigma(x,t_1) \, dx\\
&-\frac12\int_{t_1}^{t_2}\int_{\mathbb R^3} |B|^2\partial_t\sigma \, dxdt;
\end{split}
\end{equation}
while the second integral in \eqref{en-local1} can be written as 
\begin{equation}\notag
\begin{split}
\int_{t_1}^{t_2}\int_{\mathbb R^3}  ((\nabla\times B)\times B)\cdot \nabla\times(B\sigma) \, dxdt=&\int_{t_1}^{t_2}\int_{\mathbb R^3}  ((\nabla\times B)\times B)\cdot \nabla\times B\sigma \, dxdt\\
&-\int_{t_1}^{t_2}\int_{\mathbb R^3}  ((\nabla\times B)\times B)\cdot (B\times\nabla\sigma) \, dxdt\\
=&-\int_{t_1}^{t_2}\int_{\mathbb R^3}  ((\nabla\times B)\times B)\cdot (B\times\nabla\sigma) \, dxdt\\
\end{split}
\end{equation}
since the first integral on the right hand side vanishes. Therefore the local energy equality \eqref{en-local1} turns into
\begin{equation}\label{en-local2}
\begin{split}
&\int_{\mathbb R^3} |B(x,t_2)|^2\sigma(x,t_2) \, dx-\int_{\mathbb R^3} |B(x,t_1)|^2\sigma(x,t_1) \, dx\\
=&\int_{t_1}^{t_2}\int_{\mathbb R^3} |B|^2\partial_t\sigma \, dxdt+\int_{t_1}^{t_2}\int_{\mathbb R^3}  ((\nabla\times B)\times B)\cdot (B\times\nabla\sigma) \, dxdt.
\end{split}
\end{equation}
In \eqref{en-local2}, we choose $\sigma$ being radial such that $\sigma\geq 0$ and 
\begin{equation}\notag
\sigma(r)=
\begin{cases}
1, \ \ 0\leq r\leq \frac12,\\
0, \ \ r>1.
\end{cases}
\end{equation}
In particular, $\partial_t \sigma=0$. In term of self-similar form and self-similar variable $y=\frac{x}{(T-t)^{\frac1{\alpha+2}}}$, the local energy equality \eqref{en-local2} becomes
\begin{equation}\label{en-local3}
\begin{split}
&\ t_2^{\frac{3-2\alpha}{\alpha+2}} \int_{|y|\leq t_2^{-\frac1{\alpha+2}}} |H(y)|^2\sigma(yt_2^{\frac1{\alpha+2}})\, dy\\
=&\ t_1^{\frac{3-2\alpha}{\alpha+2}} \int_{|y|\leq t_1^{-\frac1{\alpha+2}}} |H(y)|^2\sigma(yt_1^{\frac1{\alpha+2}})\, dy\\
&+\int_{t_1}^{t_2}t^{\frac{2-3\alpha}{\alpha+2}}\int_{\mathbb R^3}((\nabla\times H)\times H)\cdot \left(H\times\nabla\sigma(yt^{\frac1{\alpha+2}}) \right) \, dydt.
\end{split}
\end{equation}
In view of the fact $\nabla\sigma(r)=0$ for $r<\frac12$ and $r>1$, by changing the order of the integrals, the last term of \eqref{en-local3} can be written as
\begin{equation}\notag
\int_{\frac12t_2^{-\frac1{\alpha+2}}\leq |y|\leq t_1^{-\frac1{\alpha+2}}} ((\nabla\times H)\times H)\cdot H\times \int_{t_1}^{t_2} t^{\frac{2-3\alpha}{\alpha+2}}\nabla \sigma(yt^{\frac1{\alpha+2}})\, dt dy
\end{equation}
where we have the estimate for the inner integral
\begin{equation}\notag
 \int_{t_1}^{t_2} t^{\frac{2-3\alpha}{\alpha+2}}|\nabla \sigma(yt^{\frac1{\alpha+2}})|\, dt\lesssim \int_{\frac12|y|^{-\alpha-2}}^{|y|^{-\alpha-2}}  t^{\frac{2-3\alpha}{\alpha+2}}\, dt\lesssim |y|^{-4+2\alpha}.
\end{equation}
Hence it follows from \eqref{en-local3} that
\begin{equation}\label{en-local4}
\begin{split}
&\left| t_2^{\frac{3-2\alpha}{\alpha+2}} \int_{|y|\leq t_2^{-\frac1{\alpha+2}}} |H(y)|^2\sigma(yt_2^{\frac1{\alpha+2}})\, dy-t_1^{\frac{3-2\alpha}{\alpha+2}} \int_{|y|\leq t_1^{-\frac1{\alpha+2}}} |H(y)|^2\sigma(yt_1^{\frac1{\alpha+2}})\, dy\right|\\
\leq&\ C \int_{\frac12t_2^{-\frac1{\alpha+2}}\leq |y|\leq t_1^{-\frac1{\alpha+2}}} \frac{|\nabla\times H| |H|^2}{|y|^{4-2\alpha}}\,dy.
\end{split}
\end{equation}
Taking $\ell_1= t_2^{-\frac1{\alpha+2}}$ and $\ell_2= t_1^{-\frac1{\alpha+2}}$ (with $\ell_1<\ell_2$) in \eqref{en-local4} we have
\begin{equation}\label{en-local5}
\begin{split}
&\left| \frac1{\ell_1^{3-2\alpha}} \int_{|y|\leq \ell_1} |H(y)|^2\sigma(\frac{y}{\ell_1})\, dy- \frac1{\ell_2^{3-2\alpha}} \int_{|y|\leq \ell_2} |H(y)|^2\sigma(\frac{y}{\ell_2})\, dy\, dy\right|\\
\leq&\ C \int_{\frac12\ell_1\leq |y|\leq \ell_2} \frac{|\nabla\times H| |H|^2}{|y|^{4-2\alpha}}\,dy.
\end{split}
\end{equation}

%\subsection{Case of $\alpha>\frac32$}
%It is known that $B\in L_t^\infty L_x^2$. In term of the self-similar profile we have
%\begin{equation}\notag
%\begin{split}
%\int_{\mathbb R^3}|B(t)|^2\,dx=&\ (T-t)^{-\frac{2\alpha}{\alpha+2}} \int_{|y|\leq \rho_0 (T-t)^{-\frac1{\alpha+2}}} H^2(y) (T-t)^{\frac3{\alpha+2}}\, dy\\
%=&\ L^{2\alpha-3}  \int_{|y|\leq \rho_0 L} H^2(y)\, dy
%\end{split}
%\end{equation}
%with $L=(T-t)^{-\frac1{\alpha+2}}$. Thus,
%\begin{equation}\label{L2}
% \int_{|y|\leq \rho_0 L} H^2(y)\, dy\leq C L^{3-2\alpha}.
%\end{equation}
%If $\alpha>\frac32$, taking $L\to \infty$ in \eqref{L2}, it implies $H\equiv 0$.

\medskip

\subsection{Case of $-2<\alpha\leq \frac3p$}

\begin{Lemma}\label{le-vanish}
Let $p>2$ and $-2<\alpha\leq \frac3p$. Assume $H\in L^p\cap C_{loc}^2$. We have
\[\frac{1}{\ell_2^{3-2\alpha}} \int_{|y|\leq \ell_2}|H(y)|^2 \sigma(\frac{y}{\ell_2})\, dy\to 0 \ \ \ \mbox{as} \ \ \ell_2\to\infty.\]
\end{Lemma}
\pf
For $1<M<\ell_2$ we rewrite the integral
\begin{equation}\notag
\begin{split}
&\frac{1}{\ell_2^{3-2\alpha}} \int_{|y|\leq \ell_2}|H(y)|^2 \sigma(\frac{y}{\ell_2})\, dy\\
=&\frac{1}{\ell_2^{3-2\alpha}} \int_{|y|\leq M}|H(y)|^2 \sigma(\frac{y}{\ell_2})\, dy+\frac{1}{\ell_2^{3-2\alpha}} \int_{M<|y|\leq \ell_2}|H(y)|^2 \sigma(\frac{y}{\ell_2})\, dy\\
\leq &\frac{1}{\ell_2^{3-2\alpha}} \int_{|y|\leq M}|H(y)|^2 \, dy+\frac{1}{\ell_2^{3-2\alpha}} \int_{M<|y|\leq \ell_2}|H(y)|^2\, dy.
\end{split}
\end{equation}
Since $\alpha\leq \frac3p<\frac32$, letting $\ell_2\to \infty$ yields 
\[\frac{1}{\ell_2^{3-2\alpha}} \int_{|y|\leq M}|H(y)|^2 \, dy\to 0.\]
On the other hand, applying H\"older's inequality we have
\begin{equation}\notag
\frac{1}{\ell_2^{3-2\alpha}} \int_{M<|y|\leq \ell_2}|H(y)|^2\, dy\lesssim \ell_2^{2\alpha-\frac6p} \left( \int_{M<|y|} |H(y)|^p\, dy \right)^{\frac2p} \to 0
\end{equation}
as $M\to \infty$, since $H\in L^p$ and $2\alpha-\frac6p\leq 0$. 

\cbdu

Recalling the fact $\sigma(r)=1$ on $r<\frac12$, applying Lemma \ref{le-vanish} to the local energy inequality \eqref{en-local5} we obtain
\begin{equation}\label{L2-iter}
\frac1{L^{3-2\alpha}} \int_{|y|\leq L} |H(y)|^2\, dy\leq C \int_{|y|\geq L} \frac{|\nabla\times H| |H|^2}{|y|^{4-2\alpha}}\, dy, \ \ \mbox{with} \ \ L=\frac12\ell_1.
\end{equation}

Assume $\nabla\times H\in L^q$ for $q>1$ and $H\in L^p$ for $p>2$. Applying H\"older's inequality we deduce
\begin{equation}\notag
\begin{split}
 \int_{|y|\geq L} \frac{|\nabla\times H| |H|^2}{|y|^{4-2\alpha}}\, dy\lesssim & \left(\int_{|y|\geq L}|\nabla\times H|^q\, dy\right)^{\frac1q}\left(\int_{|y|\geq L}|H|^p\, dy\right)^{\frac2p}\\
 &\cdot\left(\int_{|y|\geq L}|y|^{(2\alpha-4)\frac{pq}{pq-p-2q}}\, dy\right)^{1-\frac1q-\frac2p}\\
 \lesssim& \left(\int_L^\infty r^{(2\alpha-4)\frac{pq}{pq-p-2q}} r^2\, dr\right)^{1-\frac1q-\frac2p}\\
 \lesssim& \ L^{2\alpha-1-\frac3q-\frac6p}
\end{split}
\end{equation}
where in the last step we used the fact $2\alpha-1<\frac3q+\frac6p$. It then follows from \eqref{L2-iter} that
\begin{equation}\label{L2-base}
 \int_{|y|\leq L} |H(y)|^2\, dy\leq C L^{2-\frac3q-\frac6p}.
\end{equation}
Denote $\beta_{p,q}=2-\frac3q-\frac6p$. If $\beta_{p,q}<0$, taking $L\to\infty$ in \eqref{L2-base} yields $H\equiv 0$. 

If $\beta_{p,q}\geq 0$, we note $2<\frac{2q}{q-1}<p$. Interpolating between $L^2$ and $L^p$ we obtain an estimate for the $L^{\frac{2q}{q-1}}$ norm 
\begin{equation}\notag
\begin{split}
\int_{|y|\leq L} |H(y)|^{\frac{2q}{q-1}}\, dy=& \int_{|y|\leq L} |H(y)|^a |H(y)|^{\frac{2q}{q-1}-a}\, dy\\
\leq&\ C \left(\int_{|y|\leq L}|H(y)|^2\, dy \right)^{\frac{a}2} \left(\int_{|y|\leq L}|H(y)|^p\, dy \right)^{\frac{\frac{2q}{q-1}-a}p}
\end{split}
\end{equation}
with $a=\frac{2p}{p-2}-\frac{4q}{(q-1)(p-2)}$ (such that $\frac{a}2+ \frac{\frac{2q}{q-1}-a}p=1$). In view of \eqref{L2-base} we have 
\begin{equation}\label{inter-base}
\left( \int_{|y|\leq L} |H(y)|^{\frac{2q}{q-1}}\, dy\right)^{\frac{q-1}{q}}\leq C L^{\beta_{p,q}\alpha_{p,q}}
\end{equation}
with $\alpha_{p,q}=\frac{a(q-1)}{2q}=1-\frac p{q(p-2)}<1$.

Rearranging the right hand side of \eqref{L2-iter} and applying H\"older's inequality and \eqref{inter-base} we infer
\begin{equation}\notag
\begin{split}
\frac1{L^{3-2\alpha}}\int_{|y|\leq L} H^2(y)\, dy\leq& \frac{C}{L^{4-2\alpha}}\sum_{k=0}^\infty \frac1{2^{k(4-2\alpha)}} \int_{2^kL\leq |y|\leq 2^{k+1}L}|\nabla\times H| |H|^2 \, dy\\
\leq&  \frac{C}{L^{4-2\alpha}}\sum_{k=0}^\infty \frac1{2^{k(4-2\alpha)}} \left(\int_{2^kL\leq |y|\leq 2^{k+1}L}|\nabla\times H|^q \, dy\right)^{\frac1q}\\
&\cdot \left(\int_{2^kL\leq |y|\leq 2^{k+1}L}|H|^{\frac{2q}{q-1}} \, dy\right)^{\frac{q-1}q}\\
\leq&  \frac{C}{L^{4-2\alpha}}\sum_{k=0}^\infty \frac1{2^{k(4-2\alpha)}}\cdot (2^kL)^{\beta_{p,q}\alpha_{p,q}}\\
\leq& CL^{\beta_{p,q}\alpha_{p,q}-4+2\alpha} \sum_{k=0}^\infty 2^{k(\beta_{p,q}\alpha_{p,q}-4+2\alpha)}\\
\leq& CL^{\beta_{p,q}\alpha_{p,q}-4+2\alpha}
\end{split}
\end{equation}
where we verified that
\begin{equation}\notag
\begin{split}
\beta_{p,q}\alpha_{p,q}-4+2\alpha=&(2-\frac3q-\frac6p)(1-\frac p{q(p-2)})-4+2\alpha\\
<&\ 2-\frac3q-\frac6p-4+2\alpha\\
<&\ 2-\frac3q-\frac6p-4+\frac6p\\
<&\ 0
\end{split}
\end{equation}
since $\alpha\leq \frac3p$ and $q>1$. It then follows 
\begin{equation}\label{L2-improve}
 \int_{|y|\leq L} |H(y)|^2\, dy\leq C L^{\beta_{p,q}\alpha_{p,q}-1}.
\end{equation}
It is obvious that $\beta_{p,q}\alpha_{p,q}-1<\beta{p,q}$. Hence compared to \eqref{L2-base}, the $L^2$ estimate in \eqref{L2-improve} is improved. Consequently, the $L^{\frac{2q}{q-1}}$ estimate in \eqref{inter-base} can be improved to 
\begin{equation}\label{inter-improve}
\begin{split}
\left( \int_{|y|\leq L} |H(y)|^{\frac{2q}{q-1}}\, dy\right)^{\frac{q-1}{q}}\leq&\ C\left(\int_{|y|\leq L} |H(y)|^2\, dy\right)^{\frac{a(q-1)}{2q}}\\
\leq&\ C L^{\alpha_{p,q}(\beta_{p,q}\alpha_{p,q}-1)}.
\end{split}
\end{equation}
Applying the improved estimates \eqref{L2-improve} and \eqref{inter-improve} to the local energy inequality as before yields
\begin{equation}\notag
 \int_{|y|\leq L} |H(y)|^2\, dy\leq C L^{\beta_{p,q}\alpha_{p,q}^2-\alpha_{p,q}-1}
\end{equation}
which in turn leads to 
\begin{equation}\notag
\left( \int_{|y|\leq L} |H(y)|^{\frac{2q}{q-1}}\, dy\right)^{\frac{q-1}{q}}
\leq C L^{\alpha_{p,q}(\beta_{p,q}\alpha_{p,q}^2-\alpha_{p,q}-1)}.
\end{equation}
Repeating $n$ times of the process of applying improved estimates of $L^2$ and $L^{\frac{2q}{q-1}}$ to \eqref{L2-iter} we obtain 
\begin{equation}\label{L2-n}
 \int_{|y|\leq L} |H(y)|^2\, dy\leq C L^{\beta_{p,q}\alpha_{p,q}^n-\alpha_{p,q}^{n-1}-\cdot\cdot\cdot-1}.
\end{equation}
Since $0<\alpha_{p,q}<1$, it is obvious that 
\[\beta_{p,q}\alpha_{p,q}^n-\alpha_{p,q}^{n-1}-\cdot\cdot\cdot-1<0\]
for large enough $n$. Therefore it follows from \eqref{L2-n} that $H\equiv 0$ by taking the limit $L\to \infty$.

\subsection{Case of $\frac32<\alpha<\infty$}
It follows from \eqref{en-local5}
\begin{equation}\notag
\begin{split}
&\frac1{\ell_2^{3-2\alpha}} \int_{|y|\leq \ell_2} |H(y)|^2\, dy\\
\lesssim& \frac1{\ell_1^{3-2\alpha}} \int_{|y|\leq \ell_1} |H(y)|^2\, dy+ C \int_{\frac12\ell_1\leq |y|\leq \ell_2} \frac{|\nabla\times H| |H|^2}{|y|^{4-2\alpha}}\,dy.
\end{split}
\end{equation}
Taking $\ell_2=2L\gg 2$ and $\ell_1=2$ in the inequality above gives
\begin{equation}\label{en-local8}
\int_{|y|\leq L} |H(y)|^2\, dy
\lesssim L^{3-2\alpha}+ L^{3-2\alpha} \int_{1\leq |y|\leq 2L} \frac{|\nabla\times H| |H|^2}{|y|^{4-2\alpha}}\,dy.
\end{equation}
Applying H\"older's inequality again to the last integral we have
\begin{equation}\notag
\begin{split}
 \int_{1\leq |y|\leq 2L} \frac{|\nabla\times H| |H|^2}{|y|^{4-2\alpha}}\, dy\lesssim & \left(\int_{1\leq |y|\leq 2L}|\nabla\times H|^q\, dy\right)^{\frac1q}\left(\int_{1\leq |y|\leq 2L}|H|^p\, dy\right)^{\frac2p}\\
 &\cdot\left(\int_{1\leq |y|\leq 2L}|y|^{(2\alpha-4)\frac{pq}{pq-p-2q}}\, dy\right)^{1-\frac1q-\frac2p}\\
 \lesssim& \left(\int_1^{2L} r^{(2\alpha-4)\frac{pq}{pq-p-2q}} r^2\, dr\right)^{1-\frac1q-\frac2p}\\
 \lesssim& \max\{1, L^{2\alpha-1-\frac3q-\frac6p}\}.
\end{split}
\end{equation}
If $2\alpha-1-\frac3q-\frac6p\leq 0$, we have 
\begin{equation}\notag
\int_{|y|\leq L} |H(y)|^2\, dy \lesssim L^{3-2\alpha}
\end{equation}
which implies $H\equiv 0$ since $\alpha>\frac32$. Otherwise, if $2\alpha-1-\frac3q-\frac6p> 0$, we have
\begin{equation}\notag
\int_{|y|\leq L} |H(y)|^2\, dy \lesssim L^{3-2\alpha}+L^{2-\frac3q-\frac6p} \lesssim L^{\beta{p,q}}.
\end{equation}
As a consequence, similarly as before, we obtain
\begin{equation}\notag
\left( \int_{|y|\leq L} |H(y)|^{\frac{2q}{q-1}}\, dy\right)^{\frac{q-1}{q}}\lesssim L^{\beta_{p,q}\alpha_{p,q}}.
\end{equation}
Applying the last estimate above to \eqref{en-local8} yields
\begin{equation}\label{L2-improve-greater}
\begin{split}
&\int_{|y|\leq L} |H(y)|^2\, dy\\
\lesssim&\ L^{3-2\alpha}+ L^{-1}\sum_{k=-1}^{[\log_2L]}2^{k(4-2\alpha)} \int_{2^{-(k+1)}L\leq |y|\leq 2^{-k}L} |\nabla\times H| |H|^2\,dy\\
\lesssim&\ L^{3-2\alpha}+ L^{-1}\sum_{k=-1}^{[\log_2L]}2^{k(4-2\alpha)} \left(\int_{2^{-(k+1)}L\leq |y|\leq 2^{-k}L} |H|^{\frac{2q}{q-1}}\,dy\right)^{\frac{q}{q-1}}\\
\lesssim&\ L^{3-2\alpha}+ L^{\beta_{p,q}\alpha_{p,q}-1}\sum_{k=-1}^{[\log_2L]}2^{k(4-2\alpha-\beta_{p,q}\alpha_{p,q})}.
\end{split}
\end{equation}
If $4-2\alpha-\beta_{p,q}\alpha_{p,q}\geq 0$, it follows from \eqref{L2-improve-greater}
\begin{equation}\notag
\int_{|y|\leq L} |H(y)|^2\, dy\lesssim L^{3-2\alpha}+L^{3-2\alpha} \log_2L \to 0, \ \ \mbox{as} \ \ L\to \infty,
\end{equation}
which indicates $H\equiv 0$. 
On the other hand, if $4-2\alpha-\beta_{p,q}\alpha_{p,q}< 0$, we deduce from \eqref{L2-improve-greater}
\begin{equation}\notag
\int_{|y|\leq L} |H(y)|^2\, dy\lesssim L^{3-2\alpha}+L^{\beta_{p,q}\alpha_{p,q}-1}.
\end{equation}
If $\beta_{p,q}\alpha_{p,q}-1<0$, the proof is done; otherwise, we have 
\begin{equation}\notag
\int_{|y|\leq L} |H(y)|^2\, dy\lesssim L^{\beta_{p,q}\alpha_{p,q}-1}.
\end{equation}
Iterating the process above $n$ times gives 
\begin{equation}\notag
\int_{|y|\leq L} |H(y)|^2\, dy\lesssim L^{3-2\alpha}+L^{\beta_{p,q}\alpha_{p,q}^n-\alpha_{p,q}^{n-1}-\cdot\cdot\cdot-1}
\end{equation}
until $\beta_{p,q}\alpha_{p,q}^n-\alpha_{p,q}^{n-1}-\cdot\cdot\cdot-1<0$. Then we conclude the proof of Theorem \ref{thm-less32}.

\medskip

\subsection{Case of $\alpha=\frac32$}
Taking $\alpha=\frac32$, $\ell_1=L$ and $\ell_2=4L$ in \eqref{en-local5}, the energy inequality becomes 
\begin{equation}\label{L2-base-2}
\int_{L\leq |y|\leq 2L}H^2(y)\, dy\leq C\int_{\frac12L\leq |y|\leq 4L} \frac{|\nabla\times H||H|^2}{|y|}\, dy.
\end{equation}
It follows from \eqref{L2-base-2} and \eqref{assu4} that
\begin{equation}\label{L2-imp-2}
\begin{split}
\int_{L\leq |y|\leq 2L}H^2(y)\, dy\leq&\ C\int_{\frac12L\leq |y|\leq 4L}|y|^{-\delta}H^2(y)\, dy\\
\leq&\ C|L|^{-\delta}\int_{\frac12L\leq |y|\leq 4L}H^2(y)\, dy\\
\leq&\ C|L|^{-\delta}\left( \int_{\frac12L\leq |y|\leq L}H^2(y)\, dy+\int_{L\leq |y|\leq 2L}H^2(y)\, dy\right.\\
&\left.+\int_{2L\leq |y|\leq 4L}H^2(y)\, dy \right).
\end{split}
\end{equation}
Applying the estimate above on the interval $[\frac12L, L]$ gives
\begin{equation}\notag
\begin{split}
\int_{\frac12L\leq |y|\leq L}H^2(y)\, dy\leq&\ C|L|^{-\delta}\left( \int_{\frac14L\leq |y|\leq \frac12L}H^2(y)\, dy+\int_{\frac12L\leq |y|\leq L}H^2(y)\, dy\right.\\
&\left.+\int_{L\leq |y|\leq 2L}H^2(y)\, dy \right).
\end{split}
\end{equation}
Similar estimates can be obtained on $[L, 2L]$ and $[2L, 4L]$ for the other two integrals on the right hand side of \eqref{L2-imp-2}. After iterating this process $n$ times we obtain
\begin{equation}\notag
\begin{split}
\int_{\frac12L\leq |y|\leq L}H^2(y)\, dy\leq C_n|L|^{-n\delta}\sum_{k_1,..., k_n=-1}^{-1}\int_{2^{k_1+\cdot\cdot\cdot+k_n}L\leq |y|\leq 2^{k_1+\cdot\cdot\cdot+k_n+1}L}H^2(y)\, dy
\end{split}
\end{equation}
with a constant $C_n$ depending on $n$. Hence for any $N\in \mathbb N$, there exists a constant $C_N$ depending on $N$ such that
\begin{equation}\label{L2-final}
\int_{L\leq |y|\leq 2L}H^2(y)\, dy\leq C_NL^{-N}, \ \ \mbox{for} \ \ L\gg 1.
\end{equation}
On the other hand, thanks to \eqref{assu3}, we have
\begin{equation}\label{L2-final2}
\int_{L\leq |y|\leq 2L}H^2(y)\, dy\geq c\int_{L\leq |y|\leq 2L} |y|^{-3-2\delta}\, dy\geq c L^{-2\delta}, \ \ \mbox{for} \ \ L\gg 1.
\end{equation}
It is obvious that \eqref{L2-final2} contradicts \eqref{L2-final}. Hence we conclude $H\equiv 0$.

\bigskip

\section{Pseudo-selfsimilar solutions}
\label{sec-pseudo}
We prove Theorem \ref{thm-pse1} and Theorem \ref{thm-pse2} in this section.
Consider the solution of the electron MHD in the pseudo-selfsimilar form
\[B_{\lambda,\mu}(x,t)=\mu(t)H(\lambda(t)x)\]
with $H$ satisfying \eqref{eq-H2}. 
Taking inner product of \eqref{eq-H2} with $H$ we obtain
\begin{equation}\notag
\frac{\mu'}{\lambda^2\mu^2}\int_{\mathbb R^3}|H|^2\, dy+\frac{\lambda'}{\lambda^3\mu}\int_{\mathbb R^3}y\cdot \nabla \frac{|H|^2}2\, dy+\int_{\mathbb R^3}\nabla\times ((\nabla\times H)\times H)\cdot H\, dy=0.
\end{equation}
Applying integration by parts to the second and third integrals on the left hand side of the above equation yields
\begin{equation}\label{cond-0}
\left(\frac{\mu'}{\lambda^2\mu^2}- \frac{3\lambda'}{2\lambda^3\mu}\right) \int_{\mathbb R^3}|H|^2\, dy=0.
\end{equation}

In view of \eqref{eq-H2}, we must have that $\frac{\mu'}{\lambda^2\mu^2}H+\frac{\lambda'}{\lambda^3\mu}y\cdot\nabla H$ is independent of the time variable, which can be ensured under one of the two options:
\begin{itemize}
\item [(I)]
\begin{equation}\label{coeff}
\frac{\mu'}{\lambda^2\mu^2}=C_1, \ \ \frac{\lambda'}{\lambda^3\mu}=C_2
\end{equation}
for some constants $C_1$ and $C_2$; \\
\item [(II)] the coefficient in \eqref{cond-0} vanishes, i.e. 
\begin{equation}\label{coeff2}
\frac{\mu'}{\lambda^2\mu^2}-\frac{3\lambda'}{2\lambda^3\mu}=0.
\end{equation}
\end{itemize}

Since $\nabla^2B(x,t)=\lambda^2(t)\mu(t) H(\lambda(t)x)$, the assumption (\ref{assu2}) implies
\begin{equation}\label{cond-1}
\lim_{t\to T_0^-} \lambda^2(t)\mu(t) =\infty.
\end{equation}

In case of (I), it follows from \eqref{coeff} that
\[\frac{\mu'}{\mu}=C_1\lambda^2\mu=\frac{C_1}{C_2}\cdot\frac{\lambda'}{\lambda}.\]
Consequently we have
\begin{equation}\label{cond-2}
\mu(t)=C_0\lambda^\beta(t)
\end{equation}
for some constant $C_0>0$, with $\beta=\frac{C_1}{C_2}$.
Since
\begin{equation}\notag
\|B(t)\|_{L^2(\mathbb R^3)}=\mu(t)\lambda^{-\frac32}(t)\|H\|_{L^2(\mathbb R^3)},
\end{equation}
 assumption \eqref{assu1} implies
 \begin{equation}\label{cond-3}
 \lim_{t\to T^-} \mu(t)\lambda^{-\frac32}(t)<\infty.
 \end{equation}
 Combining \eqref{cond-1}, \eqref{cond-2} and \eqref{cond-3} we have
 \[  \lim_{t\to T^-}\lambda(t)=\infty\]
 and $-2<\beta\leq \frac32$.
 
 On the other hand, in case of (II), it follows from \eqref{coeff2} that
 \[\frac{\mu'}{\mu}=\frac32\cdot \frac{\lambda'}{\lambda}\]
 and hence 
 \[\mu(t)=C_0\lambda^{\frac32}(t)\]
 for a constant $C_0>0$. 
 
 Therefore we proceed in different cases: 
 \begin{itemize}
 \item [(i)] $\beta\in[0,\frac32)$;
 \item [(ii)] $\beta\in(-2,0)$;
 \item [(iii)] $\beta=\frac32$.
 \end{itemize}
 
{\textbf{ Case (i): $\beta=\frac{C_1}{C_2}\in [0,\frac32)$. }} Without loss of generality, we assume $C_1\geq 0$ and $C_2>0$. Recall 
\[\lambda'=C_2\lambda^3\mu=C_0C_2\lambda^{3+\beta}\]
following which we have
\begin{equation}\label{sol-la}
-(\beta+2)\lambda^{-\beta-2}(t)+(\beta+2)\lambda^{-\beta-2}(t_0)=C_0C_2(t-t_0), \ \ 0\leq t_0<t<T.
\end{equation}
Taking the limit $t\to T$ in \eqref{sol-la} and using \eqref{cond-1} yields
\begin{equation}\label{cond-4}
(\beta+2)\lambda^{-\beta-2}(t_0)=C_0C_2(T-t_0).
\end{equation}
Combining \eqref{sol-la} and \eqref{cond-4} we get
\[\lambda(t)=[(\beta+2)^{-1}C_0C_2(T-t)]^{-\frac{1}{\beta+2}}\]
and hence
\[\mu(t)=C_0\lambda^{\beta}(t)=C_0[(\beta+2)^{-1}C_0C_2(T-t)]^{-\frac{\beta}{\beta+2}}.\]
It is then easy to verify that the coefficient in \eqref{cond-0} satisfies
\[\frac{\mu'}{\lambda^2\mu^2}-\frac{3\lambda'}{2\lambda^3\mu}=(\beta-\frac32)C_0^{-1}\lambda^{-\beta-3}(t)\lambda'(t)\neq 0. \]
Consequently we infer from \eqref{cond-0} that
\[\int_{\mathbb R^3}|H|^2\, dy=0\]
and thus $H\equiv 0$. 

{\textbf{ Case (ii): $\beta=\frac{C_1}{C_2}\in (-2,0)$. }} In this case, $C_1=\frac{\mu'}{\lambda^2\mu^2}$ and $C_2=\frac{\lambda'}{\lambda^3\mu}$ have different signs. Again, it follows from \eqref{cond-0} that
$\int_{\mathbb R^3}|H|^2\, dy=0$
and $H\equiv 0$. 

{\textbf{ Case (iii): $\beta=\frac{C_1}{C_2}=\frac32$. }} We note the coefficient in \eqref{cond-0} vanishes, and hence \eqref{cond-0} is inconclusive. In this case, we make different assumptions: 

For $p>2$ and $H\not\equiv 0$, we have
\begin{equation}\notag
\|B(t)\|_{L^p(\mathbb R^3)}=\mu(t)\lambda^{-\frac3p}(t) \|H\|_{L^p(\mathbb R^3)}=C(T-t)^{-\frac37+\frac{6}{7p}}\|H\|_{L^p(\mathbb R^3)}\to \infty
\end{equation}
as $t\to T^{-}$, which contradicts the assumption \eqref{assu3}.

On the other hand, we have
\[\nabla\times J(x,t)=\mu(t)\lambda^2(t)\nabla\times \nabla \times H(y)=C(T-t)^{-1}\Delta H(y).\]
It follows, if $\Delta H\not\equiv 0$,
\begin{equation}\notag
\int_0^t \|\nabla\times J(x,\tau)\|_{L^\infty} \, d\tau=\int_0^t(T-\tau)^{-1}\, d\tau \|\Delta H\|_{L^\infty}\to \infty
\end{equation}
as $t\to T^{-}$, which contradicts the assumption \eqref{assu4}.

%\begin{Remark}
%It is not known yet whether or not the Beale-Kato-Majda blowup criterion holds for the electron MHD without resistivity. The obstacle of obtaining the BKM blowup criterion relies on the derivative loss in the nonlinear term of the electron MHD. Nevertheless, our result under assumption \eqref{assu4} provides a Beale-Kato-Majda criterion for the electron MHD without resistivity in the category of pseudo-self-similar solutions.
%\end{Remark}

%\section*{Conflict of Interest Statement}

%The author has no conflict of interests. 

%\section*{Data Availability Statement}

%My manuscript has no associated data. 

%\bigskip

%\section*{Acknowledgement}
%The author is grateful for the hospitality of the Institute for Advanced Study and Princeton University where the work was started. She would also like to thank Professor Alexandru Ionescu for inspirational conversations.

%\Endrefs
\end{document}